%% file: main.tex
\documentclass{article}

\PassOptionsToPackage{numbers, compress}{natbib}

\usepackage[final]{neurips_2024_ml4ps}

\usepackage[utf8]{inputenc} 
\usepackage[T1]{fontenc}    
\usepackage{hyperref}       
\usepackage{url}            
\usepackage{booktabs}       
\usepackage{amsmath, amsfonts}       
\usepackage{nicefrac}       
\usepackage{microtype}      
\usepackage{xcolor}         
\usepackage{multirow} 
\usepackage{graphicx}
\usepackage{caption, subcaption}
\usepackage{enumitem}
\usepackage{soul}
\usepackage{todonotes}

\setitemize{label=\textbullet, leftmargin=*, nolistsep}

\input{latex-defs}

\title{
    Scalable nonlinear manifold reduced order model \\ for dynamical systems
}

\author{%
  Ivan Zanardi${}^1$ \\
  \texttt{zanardi3@illinois.edu} \\
  \And
  Alejandro N.~Diaz${}^2$ \\
  \texttt{andiaz@sandia.gov} \\
  \And
  Seung Whan Chung${}^3$ \\
  \texttt{chung28@llnl.gov} \\
  \AND
  Marco Panesi${}^1$ \\
  \texttt{mpanesi@illinois.edu} \\
  \And
  Youngsoo Choi${}^3$ \\
  \texttt{choi15@llnl.gov} \\
  \aff
  ${}^1$University of Illinois Urbana-Champaign, Urbana, IL 61801 \\
  ${}^2$Sandia National Laboratories, Albuquerque, NM 87123 \\
  ${}^3$Lawrence Livermore National Laboratory, Livermore, CA 94550
}

\begin{document}

\maketitle

\begin{abstract}
The domain decomposition (DD) nonlinear-manifold reduced-order model (NM-ROM) represents a computationally efficient method for integrating underlying physics principles into a neural network-based, data-driven approach. Compared to linear subspace methods, NM-ROMs offer superior expressivity and enhanced reconstruction capabilities, while DD enables cost-effective, parallel training of autoencoders by partitioning the domain into algebraic subdomains. In this work, we investigate the scalability of this approach by implementing a ``bottom-up'' strategy: training NM-ROMs on smaller domains and subsequently deploying them on larger, composable ones. The application of this method to the two-dimensional time-dependent Burgers' equation shows that extrapolating from smaller to larger domains is both stable and effective. This approach achieves an accuracy of 1\% in relative error and provides a remarkable speedup of nearly 700 times.
\end{abstract}

\input{01_introduction}
\input{02_dd_fom}
\input{03_dd_nmrom}
\input{04_architecture}
\input{05_numerics}
\input{06_conclusion}

\begin{ack}
This work was performed at Lawrence Livermore National Laboratory.
\textbf{I.\ Zanardi} was supported by the Data Science Summer Institute (DSSI) and LDRD (22-SI-006) at Lawrence Livermore National Laboratory and the Vannevar Bush Faculty Fellowship OUSD(RE) Grant N00014-21-1-295.
\textbf{A.\ N.\ Diaz} was supported in part by a 2021 DoD National Defense Science and Engineering Graduate Fellowship and the S. Scott Collis Fellowship at Sandia National Laboratory.
\textbf{S.\ W.\ Chung} was supported by LDRD (22-SI-006).
\textbf{M.\ Panesi} was supported by the Vannevar Bush Faculty Fellowship OUSD(RE) Grant N00014-21-1-295.
\textbf{Y.\ Choi} was supported by the U.S. Department of Energy, Office of Science, Office of Advanced Scientific Computing Research, as part of the CHaRMNET Mathematical Multifaceted Integrated Capability Center (MMICC) program, under Award Number DE-SC0023164 and partially by LDRD (22-SI-006).
Lawrence Livermore National Laboratory is operated by Lawrence Livermore National Security, LLC, for the U.S. Department of Energy, National Nuclear Security Administration under Contract DE-AC52-07NA27344. IM release number: LLNL-CONF-869013.
\end{ack}

\bibliography{references, yc_references}

\end{document}

%% file: latex-defs.tex
\newcommand{\eref}[1]{\mbox{\rm(\ref{#1})}}

\newcommand{\set}[1]{\left\{ #1 \right\} }
\newcommand{\norm}[1]{\left\| #1 \right\|}

\newcommand{\real}{\mathbb{R}}

\newcommand {\cD} {{\cal D}}

\newcommand{\bzero}{\mathbf{0}}

\newcommand{\bff}{\mathbf{f}}
\newcommand{\bg}{\mathbf{g}}
\newcommand{\bh}{\mathbf{h}}

\newcommand{\br}{\mathbf{r}}

\newcommand{\bx}{\mathbf{x}}

\newcommand{\bA}{\mathbf{A}}
\newcommand{\bB}{\mathbf{B}}

\newcommand{\bI}{\mathbf{I}}

\newcommand{\hbx}{\widehat{\mathbf{x}}}

\newcommand{\hbA}{\widehat{\mathbf{A}}}

\newcommand{\bmu}{{\boldsymbol{\mu}}}

\newcommand{\ovdt}{\frac{d}{dt}}

%% file: 01_introduction.tex
\section{Introduction}\label{sec:introduction}
Complex tasks such as design optimization and uncertainty quantification often require repeated simulations of a large-scale, parameterized, nonlinear system, commonly referred to as the full-order model (FOM). This approach can become impractical for large-scale problems due to the computational demands involved. Model reduction addresses this challenge by substituting the FOM with a more computationally efficient, lower-dimensional model known as a reduced-order model (ROM). Despite its advantages, constructing accurate and efficient ROMs presents its own set of challenges. In this work we explore the framework proposed by Diaz \textit{et al.}~\cite{ANDiaz_YChoi_MHeinkenschloss_2024,Diaz_Thesis_2024}, which combines the nonlinear-manifold ROM (NM-ROM) approach with an algebraic domain-decomposition (DD) framework.

Various model reduction techniques have been integrated with DD, including reduced basis elements (RBE)~\cite{YMaday_EMRonquist_2002a,
YMaday_EMRonquist_2004a,
LIapichino_AQuarteroni_GRozza_2012a,
PFAntonietti_PPacciarini_AQuarteroni_2016a,
JLEftang_DBPHuynh_DJKnezevic_EMRonquist_ATPatera_2012a,
DBPHuynh_DJKnezevic_ATPatera_2013a,JLEftang_ATPatera_2013a} or the alternating Schwarz method~\cite{MBuffoni_HTelib_AIollo_2009a,
JBarnett_ITezaur_AMota_2022a,
KSmetana_TTaddei_2022a,AIollo_GSambataro_TTaddei_2022a}, which are often tailored to specific problems and address the physical domain at the PDE level. An alternative approach is the algebraic method proposed by Hoang \textit{et al.}~\cite{CHoang_YChoi_KCarlberg_2021a},
which involves decomposing the FOM at the discrete level and computing linear-subspace reduced-order models (LS-ROMs) for each subdomain. Although LS-ROMs perform effectively in many cases~\cite{BHaasdonk_2017a,
AQuarteroni_AManzoni_FNegri_2016a, 
MHinze_SVolkwein_2005a,
MGubisch_SVolkwein_2017a,
cheung2023local,
copeland2022reduced,carlberg2018conservative,
ACAntoulas_2005a,PBenner_TBreiten_2017a, 
ACAntoulas_CABeattie_SGugercin_2020a,CGu_2011a,PBenner_TBreiten_2015a,
AJMayo_ACAntoulas_2007a,
ACAntoulas_IVGosea_ACIonita_2016a,
IVGosea_ACAntoulas_2018a,
choi2021space,
kim2021efficient,
choi2019space,
cheung2023datascarce,
tsai2023accelerating,
huhn2023parametric,
mcbane2022stress,
mcbane2021component,
choi2020gradient,Chung_CMAME_2024,chung2024scaled}, they are known to struggle with advection-dominated problems and those with sharp gradients, which are characterized by slowly decaying Kolmogorov $n$-width~\cite{MOhlberger_SRave_2016a}. Recent approaches, such as nonlinear-manifold reduced-order models (NM-ROMs), address these challenges by approximating the FOM within a low-dimensional nonlinear manifold. This is typically achieved by training an autoencoder on FOM snapshot data~\cite{KKashima_2016a, DHartman_LKMestha_2017a, KLee_KTCarlberg_2020a, YKim_YChoi_DWidemann_TZohdi_2022a, kim2020efficient}. However, training NM-ROMs is computationally expensive due to the high dimensionality of the FOM training data, leading to many neural network (NN) parameters. To mitigate this, Barnett \textit{et al.}~\cite{JBarnett_CFarhat_YMaday_2023a} first computed a low-dimensional proper orthogonal decomposition (POD) model and then trained the NN on the POD coefficients. Instead, we adopt the approach by Diaz \textit{et al.}~\cite{ANDiaz_YChoi_MHeinkenschloss_2024, Diaz_Thesis_2024}, which integrates an autoencoder framework with DD. This method allows for the computation of FOM training data on subdomains, thereby reducing the dimensionality of the subdomain NM-ROM training data. Consequently, fewer parameters are required for training in each subdomain.

In this work, we adopt the DD NM-ROM approach for its proven effectiveness in tackling large-scale problems with slowly decaying Kolmogorov $n$-width. This method showed superior accuracy and performance compared to both LS-ROMs and monolithic NN-ROMs~\cite{ANDiaz_YChoi_MHeinkenschloss_2024}. We specifically examine the scalability of the framework by implementing a ``bottom-up'' training strategy. This involves using snapshots from subdomains of a small-sized domain for training, and subsequently deploying the trained autoencoders to a larger, composable domain. The architecture employed is a wide, shallow, and sparse autoencoder with a sparsity mask applied to both the encoder input layer and decoder output layer, as described in \cite{YKim_YChoi_DWidemann_TZohdi_2022a}. The DD NM-ROM approach is applied to the two-dimensional time-dependent Burgers' equation.

%% file: 02_dd_fom.tex
\section{DD full order model}\label{sec:dd_fom}
First consider the monolithic FOM as a system of ODEs
\begin{equation}\label{eq:fom}
    \ovdt \bx(t) =\bff\left(\bx(t);\bmu\right),
    \qquad t \in [0, T],
    \qquad \bx(0) =\bx_0(\bmu),
\end{equation}
where $\bx:[0, T] \rightarrow \real^{N_x}$ is the high-dimensional state, $\bmu \in \cD \subset \real^{N_u}$ is a parameter, and $\mathbf{f}: \real^{N_x} \times \cD \rightarrow \real^{N_x}$. Applying Backward Euler (BE) time stepping to \eqref{eq:fom} results in the following (nonlinear) system of equations, which must be solved at each time step $k=1, \ldots, N_t$:
\begin{equation}\label{eq:fom_monolithic}
    \br\left(\bx^{(k)}, \bx^{(k-1)}; \bmu\right) =\bx^{(k)}-\bx^{(k-1)}-\tau \mathbf{f}\left(\bx^{(k)}; \bmu\right)=\bzero,
    \qquad \bx^{(0)} =\bx_0(\bmu),
\end{equation}
where $\tau=T / N_t$ denotes the time step, and $\br: \real^{N_x} \times \real^{N_x} \times \cD \rightarrow \real^{N_x}$ is the residual function.
FOMs of the form \eref{eq:fom_monolithic} typically arise from discretizations of partial differential equations (PDEs). One can reformulate \eref{eq:fom_monolithic} into a DD formulation by partitioning the residual equation into $n_\Omega$ systems of equations (so-called {\it algebraic} subdomains), coupling them via {\it compatibility constraints}, and converting the systems of equations into a least-squares problem, resulting in
\begin{equation}\label{eq:dd_fom}
\min _{\left(\bx_i^{\Omega(k)}, \bx_i^{\Gamma(k)}\right)} \frac{h}{2} \sum_{i=1}^{n_{\Omega}}\left\|\br_i\left(\bx_i^{\Omega(k)}, \bx_i^{\Gamma(k)}, \bx_i^{\Omega(k-1)}, \bx_i^{\Gamma(k-1)};\bmu\right)\right\|_2^2,
\quad \mathrm{s.t.} \quad
\sum_{i=1}^{n_{\Omega}} \bA_i \bx_i^{\Gamma(k)}=\bzero,
\end{equation}
accompanied by the proper initial conditions $\bx_i^{\Omega(0)}$ and $\bx_i^{\Gamma(0)}$, with $h > 0$ being a scaling factor and $\bx_i^{\Omega(k)}\in \real^{N_i^\Omega}$,
$\bx_i^{\Gamma(k)}\in \real^{N_i^\Gamma}$,
$\br_i:\real^{N_i^\Omega}\times \real^{N_i^\Gamma}\times \real^{N_i^\Omega}\times \real^{N_i^\Gamma}\times \cD \to \real^{N_i^r}$, and
$\bA_i\in \set{-1, 0, 1}^{N_A\times N_i^\Gamma}$
being the $i$-th subdomain interior-state, interface-state, 
residual function, and compatibility constraint matrix, respectively. The structure of the subdomain residual functions $\br_i$ and the division of the state $\bx$ into subdomain states $\left(\bx_i^\Omega, \bx_i^\Gamma\right)$ are influenced by the sparsity pattern of the overall residual function $\br$. Specifically, the interior states $\bx_i^\Omega$ are used solely for computing the residual $\br_i$ within the $i$-th subdomain, while the interface states $\bx_i^\Gamma$ are also involved in the residual computations for adjacent subdomains. Additionally, we can define port states $\bx_j^p$ as a set of nodes uniquely shared by multiple subdomains. By combining multiple ports, we can generate the interface states $\bx_i^\Gamma=\bigcup_{j\in i}\bx_j^p$. The equality constraint determined by $\bA_i$ enforces equality on the overlapping interface states. For further details, see \cite[Sec. 2]{ANDiaz_YChoi_MHeinkenschloss_2024} or \cite[Sec. 2]{CHoang_YChoi_KCarlberg_2021a}.
\\
The dependence of the residual function on the previous time step and the parameter $\bmu$ will be omitted until required.

%% file: 03_dd_nmrom.tex
\section{DD nonlinear-manifold reduced order model}\label{sec:dd_nmrom}
In the general formulation, for each subdomain $i \in \mathcal{S}^{\Omega}=\set{1, \dots, n_\Omega}$, 
let $\bg_i^\Omega: \real^{n_i^\Omega} \to \real^{N_i^\Omega}$, $n_i^\Omega \ll N_i^\Omega$, and 
$\bg_i^\Gamma: \real^{n_i^\Gamma} \to \real^{N_i^\Gamma}$, $n_i^\Gamma \ll N_i^\Gamma$,
be decoders such that 
$\bx_i^{\Omega} \approx \bg_i^\Omega\left(\hbx_i^{\Omega}\right)$ and
$\bx_i^{\Gamma} \approx \bg_i^\Gamma\left(\hbx_i^{\Gamma}\right)$. The corresponding encoders are denoted by $\bh_i^\Omega$ and $\bh_i^\Gamma$.
Also let $\bB_i \in \set{0,1}^{N_i^B \times N_i^r}$, $N_i^B \leq N_i^r$, denote a row-sampling matrix for collocation hyper-reduction (HR).
The DD NM-ROM is evaluated by solving
\begin{equation}\label{eq:dd_nmrom}
\min_{\left(\hbx_i^{\Omega(k)}, \hbx_i^{\Gamma(k)}\right)}
\frac{h}{2}\sum_{i=1}^{n_\Omega} \norm{\bB_i \br_i\left(\bg_i^\Omega\left(\hbx_i^{\Omega(k)}\right), \bg_i^\Gamma\left(\hbx_i^{\Gamma(k)} \right)\right)}_2^2,
\quad \mathrm{s.t.} \quad \sum_{i=1}^{n_\Omega} \hbA_i \hbx_i^{\Gamma(k)} = \bzero.
\end{equation}
with $\hbx_i^{\Omega(0)} = \bh_i^\Omega\big(\bx_i^{\Omega(0)}\big)$ and $\hbx_i^{\Gamma(0)} = \bh_i^\Gamma\big(\bx_i^{\Gamma(0)}\big)$. In this work, the following set up is used: 
\begin{itemize}
    \item HR is not applied, meaning $\bB_i=\bI$. However, in most cases, applying HR is crucial for achieving greater speedups, as it avoids the need to evaluate the full residual of the FOM.
    \item The interface states encoder/decoder are formulated as a combination of port encoders/decoders, e.g., $\bg_i^\Gamma = \bigcup_{j\in i}\bg_j^p$. Similarly, the latent interface states are expressed as $\hbx_i^\Gamma=\bigcup_{j\in i}\hbx_j^p$.
    \item The compatibility constraints in \eqref{eq:dd_nmrom} are expressed at the ROM latent interface states, similar to the approach used for the FOM, with $\hbA_i \in\{-1,0,1\}^{n_A \times n_i^{\Gamma}}$. This approach is known as the strong ROM-port constraint (SRPC) formulation~\cite{ANDiaz_YChoi_MHeinkenschloss_2024}.
\end{itemize}
The DD NM-ROM formulation \eref{eq:dd_nmrom} provides several advantages. Its training process, which involves computing $\bg_i^\Omega$ and $\bg_j^p$, is localized, requires few parameters, and can be executed in parallel. This allows for the adaptation of ROMs to specific features of the problem, potentially resulting in more compact ROMs. Additionally, the framework supports parallelization to accelerate both ROM computation/training and execution, and it enables the implementation of a ``bottom-up'' training strategy, as detailed in the introductory section.
\vskip 0mm
\textbf{Solver}\hspace{2mm}
The DD FOM \eref{eq:dd_fom} and DD NM-ROM \eref{eq:dd_nmrom} are solved using an inexact Lagrange-Newton sequential quadratic programming (SQP) solver. In this approach, the Hessian of the Lagrangian is approximated using a Gauss-Newton method, thus eliminating the need for computing second-order derivatives of residuals and constraints in \eref{eq:dd_nmrom}. Despite this approximation, the method still ensures effective convergence for both \eref{eq:dd_fom} and \eref{eq:dd_nmrom}. For more details, refer to \cite{ANDiaz_YChoi_MHeinkenschloss_2024}.

%% file: 04_architecture.tex
\vskip 0mm
\textbf{Autoencoder}\hspace{2mm}
We use single-layer, wide, and sparse decoders with smooth activation functions to represent the maps $\bg_i^\Omega$ and $\bg_j^p$. 
The corresponding encoders, denoted $\bh_i^\Omega$ and $\bh_j^p$,
are also single-layer, wide, and sparse.
Shallow networks are used for computational efficiency; fewer layers correspond to fewer repeated matrix-vector multiplications when evaluating the decoders. 
The shallow depth necessitates a wide network to maintain enough expressiveness for use in NM-ROM. 
Smooth activations (i.e., \textit{swish}) are used to ensure that $\bg_i^\Omega$ and $\bg_j^p$ are continuously differentiable. Sparsity is imposed on both the input and output layers of the autoencoder to maintain symmetry across the latent layer. This sparsity pattern adopts a tri-banded structure, inspired by 2D finite difference stencils, with the number of nonzero elements per band and the spacing between bands serving as hyperparameters~\cite{ANDiaz_YChoi_MHeinkenschloss_2024}. Normalization and de-normalization layers are also applied at the encoder input and decoder output layers, respectively, and Gaussian noise is added during training.
\\
To train the autoencoders, we first gather $N_t$ snapshots of the interior and port states in an \textit{offline} stage by solving \eref{eq:fom_monolithic} for each parameter $\bmu_\ell$ for $\ell=1,\dots,M$. The snapshots are then randomly divided into an 80-20 split for training and validation. The Adam optimizer is used to minimize the MSE loss over 1000 epochs, with a batch size of 1024. We also incorporate early stopping and reduce the learning rate on plateau, starting with an initial value of $10^{-3}$. The implementation is carried out using PyTorch, leveraging the PyTorch Sparse and SparseLinear packages.

%% file: 05_numerics.tex
\section{Numerical experiment}\label{sec:numerics}
\begin{figure}[htb!]
    \centering
    \includegraphics[width=0.6\textwidth]{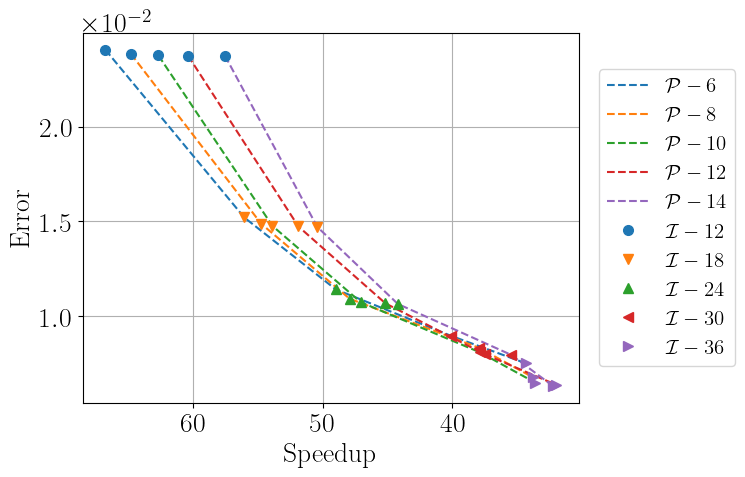}
    \caption{Averaged speedups and errors from 50 distinct $2 \times 2$ test cases, evaluated across 25 different combinations of latent dimensions for vertical/horizontal port nodes ($\mathcal{P}$) and interior nodes ($\mathcal{I}$). The results illustrate the impact of varying latent space dimensions on performance metrics.}
    \label{fig:pareto}
\end{figure}

We consider the 2D Burgers' equation
\begin{equation}\label{eq:burgers_pde_2d}
    \frac{\partial u}{\partial t} + u\frac{\partial u}{\partial x} + v \frac{\partial u}{\partial y} = \nu \left(\frac{\partial^2 u}{\partial x^2} + \frac{\partial^2 u}{\partial y^2}\right),
    \qquad
    \frac{\partial v}{\partial t} + u \frac{\partial v}{\partial x} + v \frac{\partial v}{\partial y} = \nu \left(\frac{\partial^2 v}{\partial x^2} + \frac{\partial^2 v}{\partial y^2}\right),
\end{equation}
with viscosity $\nu=10^{-3}$. The PDE is discretized using centered finite differences over a uniform structured mesh, with grid sizes $h_x$ and $h_y$. The integration is performed for $T=2$ s with a time step of $\tau = 0.02$ s. To evaluate the performances of the proposed DD NM-ROM, the following discrete $L^{\infty}-L^2$ error is used:
\begin{equation}\label{eq:relative_error}
e=\max _k\left[\frac{h_x h_y}{n_{\Omega}} \sum_{i=1}^{n_{\Omega}}\left(\left\|\mathbf{x}_i^{\Omega(k)}-\mathbf{g}_i^{\Omega}\left(\widehat{\mathbf{x}}_i^{\Omega(k)}\right)\right\|_2^2+\left\|\mathbf{x}_i^{\Gamma(k)}-\mathbf{g}_i^{\Gamma}\left(\widehat{\mathbf{x}}_i^{\Gamma(k)}\right)\right\|_2^2\right)\right]^{1 / 2} .
\end{equation}
The implementation was performed sequentially. However, to highlight the potential benefits of a parallel implementation, the reported wall clock time for computing subdomain-specific quantities in the SQP solver is based on the maximum wall clock time recorded among all subdomains. The wall clock time for the remaining steps of the SQP solver is measured as the total wall clock time.
\\
Training was conducted using the Lassen machine, while testing and computations were carried out on the Dane machine at Lawrence Livermore National Laboratory. For more details, refer to~\cite{llnl_machines}. The code can be found at https://github.com/LLNL/DD-NM-ROM.
\vskip 0mm
\textbf{Model construction}\hspace{2mm}
For training, we employ a DD problem with periodic boundary conditions, divided into four uniformly sized subdomains ($n_{\Omega} = 4$). These subdomains are arranged in a $2 \times 2$ configuration within the domain $(x, y) \in [0, 1] \times [0, 1]$, with each direction containing 100 grid nodes. We explore the effects of the scalar parameter $\mu_i$, which represents the peak value of the initial sinusoidal function in the $i$-th subdomain, as defined below:
\begin{equation}\label{eq:sampling}
\begin{aligned}
    u_i^{(0)} = v_i^{(0)} & = 
        \left| \mu_i \sin (2 \pi x) \sin (2 \pi y) \right|
    \quad \forall i \in \mathcal{S}^{\Omega}, \\
    \text{with} \quad \mu_i & = \gamma_i \xi_i,
    \quad \gamma_i \sim \mathcal{U}(0.5, 1.5),
    \quad \xi_i \sim \mathcal{B}(1, 0.5),
\end{aligned}
\end{equation}
where $\mathcal{U}$ and $\mathcal{B}$ denote uniform and binomial distributions, respectively. We sample a total of 2000 different initial $2 \times 2$ configurations from \eref{eq:sampling} using a Latin hypercube sampling (LHS) strategy, which ensures thorough exploration of all possible configurations. The sole constraint is that the bottom-left subdomain always has $\xi_0 = 1$, ensuring the model can observe waves traveling across the entire domain.
\\
In our study, we constructed three distinct autoencoders: one for the interior states $\left(\bg_i^\Omega = \bg^\Omega\right.$ $\left.\forall i \in \mathcal{S}^\Omega\right)$, and two others for the vertical $\left(\bg_\mathrm{v}^p\right)$ and horizontal $\left(\bg_\mathrm{h}^p\right)$ ports. We explored various latent space dimensions, ranging from 12 to 36 for the interior states and from 6 to 14 for both the vertical and horizontal port states. The results shown in Fig. \ref{fig:pareto} are based on averaged speedups and errors from 50 distinct $2 \times 2$ test cases, evaluated across 25 different combinations of latent dimensions for vertical/horizontal port nodes ($\mathcal{P}$) and interior nodes ($\mathcal{I}$). These results indicate that the dimensions of the latent space for the port states have a notably smaller impact on performance compared to those for the interior states, likely due to differences in compression rates.

\begin{figure}[htb!]
\centering
\begin{minipage}{\textwidth}
    \centering
    \begin{minipage}{0.07\textwidth}
        \centering \phantom{Time}
    \end{minipage}%
    \begin{minipage}{0.3\textwidth}
        \centering DD FOM
    \end{minipage}%
    \begin{minipage}{0.3\textwidth}
        \centering DD NM-ROM
    \end{minipage}%
    \begin{minipage}{0.3\textwidth}
        \centering Absolute Error
    \end{minipage}
    \\[5pt]
    \begin{minipage}{0.1\textwidth}
        \centering $t=0$ s
    \end{minipage}%
    \begin{minipage}{0.3\textwidth}
        \centering
        \includegraphics[width=0.95\textwidth]{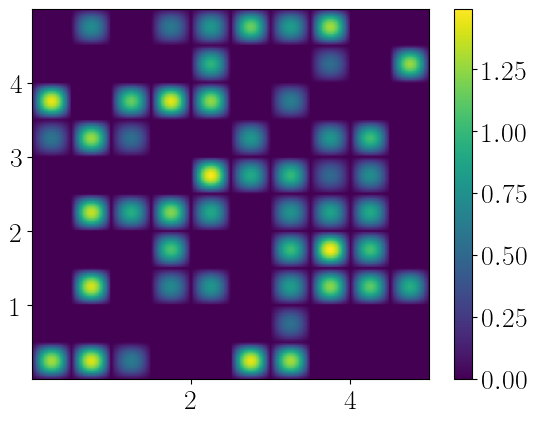}
    \end{minipage}%
    \begin{minipage}{0.3\textwidth}
        \centering
        \includegraphics[width=0.95\textwidth]{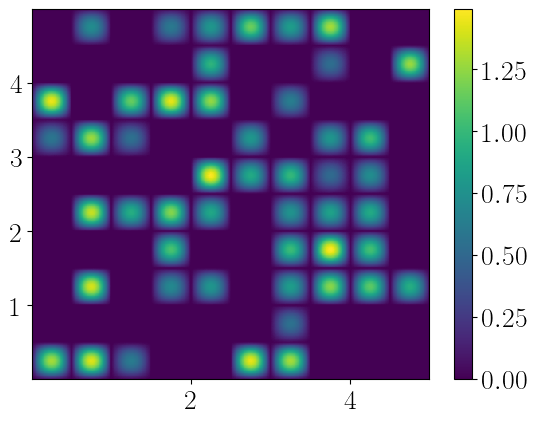}
    \end{minipage}%
    \begin{minipage}{0.3\textwidth}
        \centering
        \includegraphics[width=0.95\textwidth]{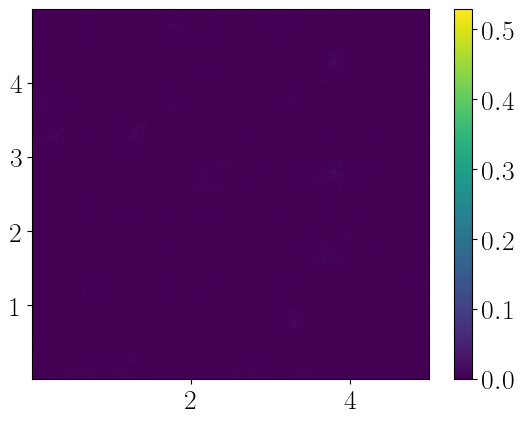}
    \end{minipage}
    \\[1pt]
    \begin{minipage}{0.1\textwidth}
        \centering $t=1$ s
    \end{minipage}%
    \begin{minipage}{0.3\textwidth}
        \centering
        \includegraphics[width=0.95\textwidth]{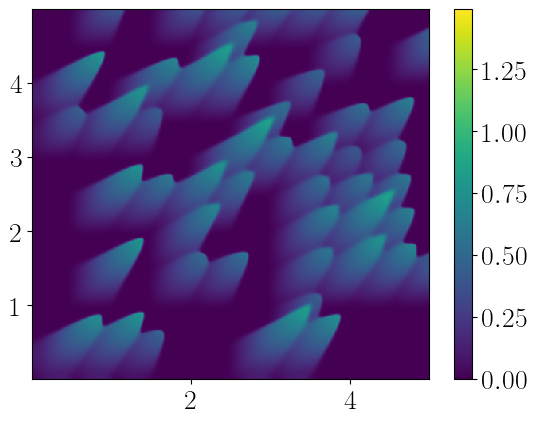}
    \end{minipage}%
    \begin{minipage}{0.3\textwidth}
        \centering
        \includegraphics[width=0.95\textwidth]{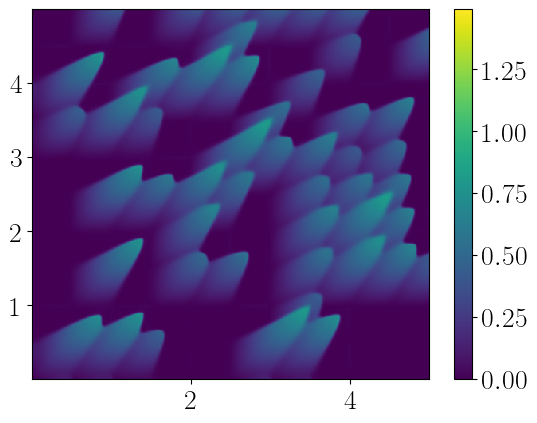}
    \end{minipage}%
    \begin{minipage}{0.3\textwidth}
        \centering
        \includegraphics[width=0.95\textwidth]{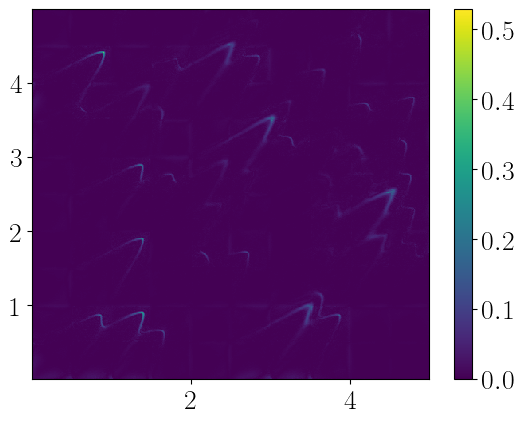}
    \end{minipage}
    \\[1pt]
    \begin{minipage}{0.1\textwidth}
        \centering $t=2$ s
    \end{minipage}%
    \begin{minipage}{0.3\textwidth}
        \centering
        \includegraphics[width=0.95\textwidth]{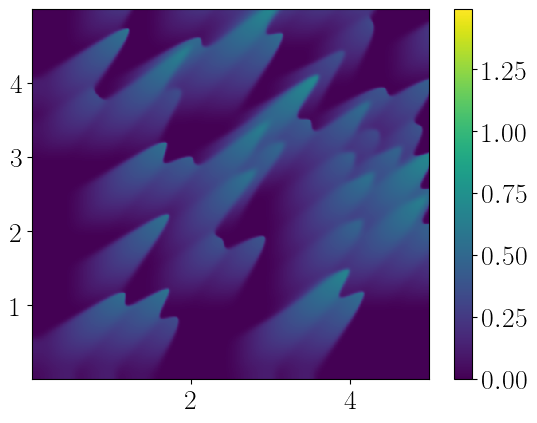}
    \end{minipage}%
    \begin{minipage}{0.3\textwidth}
        \centering
        \includegraphics[width=0.95\textwidth]{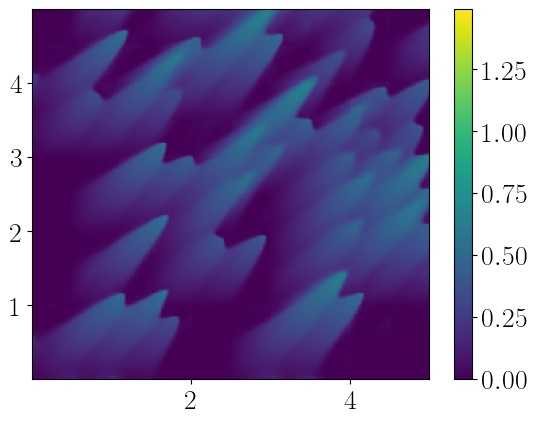}
    \end{minipage}%
    \begin{minipage}{0.3\textwidth}
        \centering
        \includegraphics[width=0.95\textwidth]{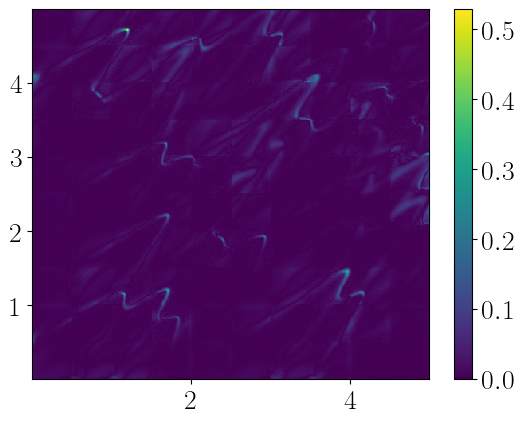}
    \end{minipage}
\end{minipage}
\caption{The $u$ velocity predicted by the DD FOM and DD NM-ROM models, along with the absolute error between them, at various time instants. The initial condition was randomly sampled using \eref{eq:sampling} on a $10 \times 10$ configuration.}
\label{fig:u_10by10}
\end{figure}

\textbf{Model deployment}\hspace{2mm}
To assess our ``bottom-up'' strategy, we applied the Pareto optimal NM-ROM from the previous section (with $\mathcal{I}=24$ and $\mathcal{P}=10$) to a DD problem with homogeneous Neumann boundary conditions and $n_{\Omega} = 100$ subdomains, each uniformly sized and arranged in a $10 \times 10$ configuration within the spatial range $(x, y) \in [0, 5] \times [0, 5]$. Each direction was discretized using 500 grid nodes and the initial condition was randomly sampled from \eref{eq:sampling}. At each time step, the decoder $\bg^\Omega$ is evaluated for each subdomain to update the latent representation of the interior states. Similarly, the decoders $\bg_\mathrm{v}^p$ and $\bg_\mathrm{h}^p$ are evaluated for each vertical and horizontal port, respectively, to update the latent representations associated with those ports. As illustrated in Fig. \ref{fig:u_10by10}, which shows the $u$ velocity component at three different time instants, the predicted solution from the DD NM-ROM approximates well the true DD FOM solution. The error, quantified at $1.21 \times 10^{-2}$, is consistent with the average errors shown in Fig. \ref{fig:pareto} for the same latent dimensions, while achieving an excellent speedup of 662.62 when compared to the DD FOM. However, increased error is observed at the shock wave fronts, attributed to accumulation error phenomena. This issue can be addressed by employing dynamic training approaches, such as adjoint-based optimization, which allows the model to consider the evolution of the approximated dynamics. The notable speedup, achieved without using HR, is attributed to the inexact SQP solver, where the size of the linear system solved at each iteration is significantly smaller for the ROM.

We want to emphasize that DD NM-ROM approach offers a significant advantage by enabling faster training on a smaller $2 \times 2$ domain, avoiding the high computational cost and time required to train a monolithic model on a larger $10 \times 10$ domain with random initial conditions.

%% file: 06_conclusion.tex
\section{Conclusion}\label{sec:conclusion}
In this study, we assessed the effectiveness of the ``bottom-up'' strategy for the DD NM-ROM framework proposed in \cite{ANDiaz_YChoi_MHeinkenschloss_2024,Diaz_Thesis_2024}. Our approach began with the development of three distinct autoencoders tailored for interior states and vertical and horizontal ports, using a small training domain with a $2\times 2$ subdomain configuration. These models were then applied to a larger composable domain featuring a $10\times 10$ subdomain configuration. The results indicate that extrapolating from the smaller to the larger domain is both stable and effective, achieving accuracy and speedup improvements, with the DD NM-ROM operating approximately 700 times faster than the DD FOM. Future research will aim to integrate adjoint-based optimization techniques to mitigate potential accumulation errors. Additionally, the framework's applicability will be extended to more complex physical systems, including the Kuramoto–Sivashinsky equation, the Korteweg–De Vries equation, and the Navier–Stokes equations.